\documentclass[a4paper,11pt]{article}

\usepackage{amsmath, amssymb, amsthm}
\usepackage{graphicx}
\usepackage[colorlinks=true, urlcolor=blue, citecolor=black, linkcolor=black]{hyperref}
\usepackage{geometry}

\title{A Minimal Mathematical Model for Conducting Patterns}
\author{Tom Verhoeff\\
  Dept.\ Math.\ \& CS, Eindhoven University of Technology, The Netherlands\\
  \href{mailto:T.Verhoeff@tue.nl}{\ttfamily T.Verhoeff@tue.nl}}
\date{08~Apr 2026}

\begin{document}

\maketitle

\begin{abstract}
We present a minimal mathematical model for conducting patterns that separates geometric trajectory from temporal parametrization.
The model is based on a cyclic sequence of preparation and ictus points connected by cubic Hermite segments with constrained horizontal tangents, combined with a quintic timing law controlling acceleration and deceleration.
A single parameter governs the balance between uniform motion and expressive emphasis.
The model provides a compact yet expressive representation of conducting gestures.
It is implemented as the interactive Wolfram Demonstration \emph{Conducting Patterns} and
is used in the \emph{Crusis} web app.
\end{abstract}

\begin{figure}[hbt]
\centering
{\includegraphics[trim=40mm 30mm 45mm 40mm,clip,height=6cm]{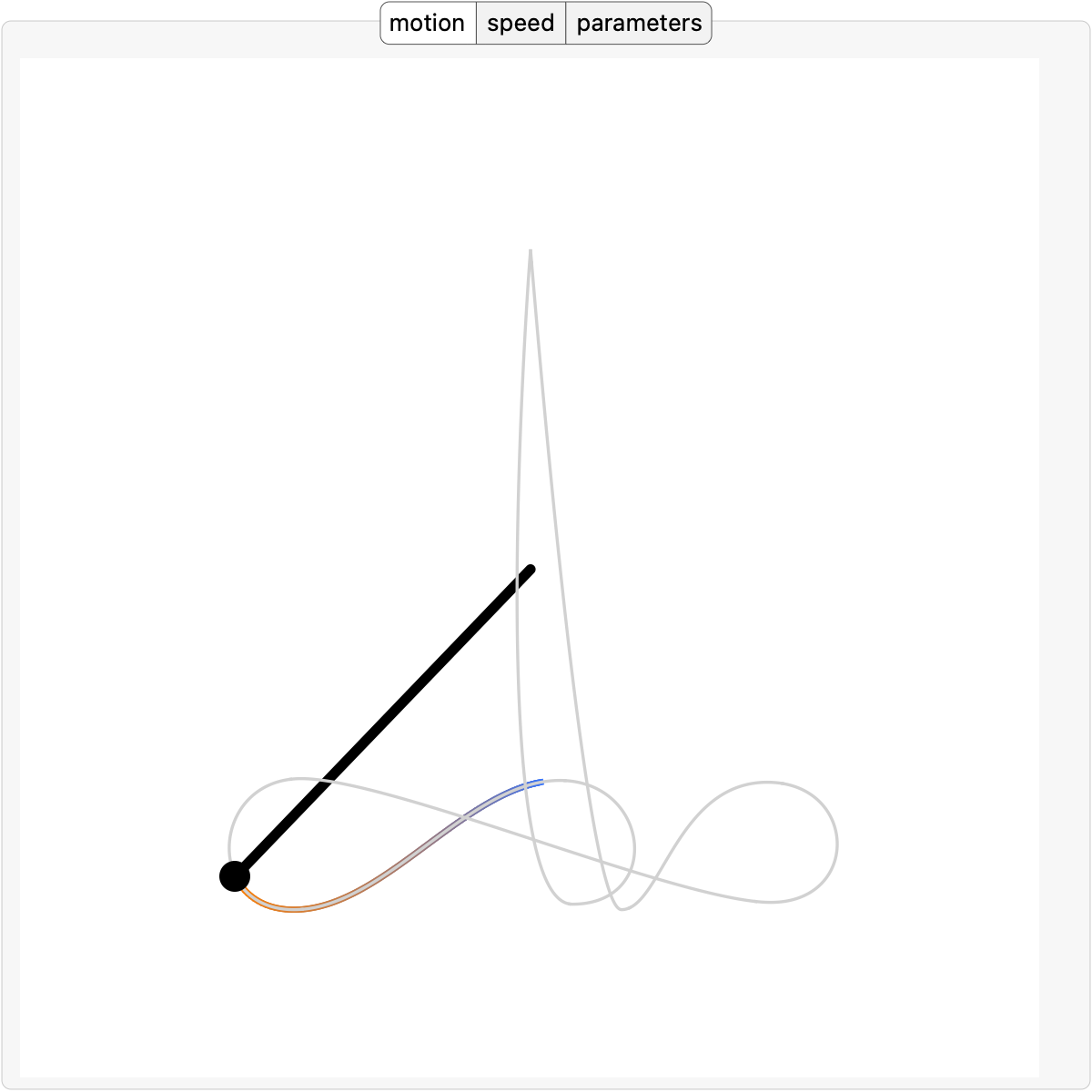}}
\caption{Baton with trail tracing a conducting pattern for a 4-beat time signature.}
\label{fig:curve+trail+baton}
\end{figure}

\section{Introduction}

Conducting gestures (see Figure~\ref{fig:curve+trail+baton}) convey temporal and structural information in music through motion.
Central to this communication are the \emph{ictus}, marking the beat,
and the \emph{preparation}, guiding performers toward it.

Despite their importance, conducting patterns are typically described pedagogically rather than mathematically.
Several digital tools attempt to support conducting practice by visualizing beat patterns or providing tempo guidance.
For example, the mobile application \emph{Maestro} (iOS/iPadOS) presents animated conducting patterns for different meters, while applications such as \emph{Pro Metronome} include visual beat indicators resembling conducting gestures.
Other educational apps provide static or animated diagrams of standard conducting patterns.
However, these tools typically rely on predefined animations rather than a parameterized mathematical model of conducting motion.

This paper proposes a minimal formal model that captures the essential structure of conducting gestures while remaining computationally simple and expressive.
Our model can be interpreted as an idealized representation of conducting gestures that abstracts away from biomechanical complexity while preserving perceptually salient features.
It aims to complement empirical and pedagogical approaches.

\section{Mathematical Model}

Our model is based on the following principles:
\begin{enumerate}
\item \textbf{Separation of concerns:} geometry and timing are modeled independently.
\item \textbf{Anchors:} each beat consists of a preparation point and an ictus.
\item \textbf{Extrema with horizontal tangents:} preparation points are local vertical maxima and ictus points local vertical minima, enforced by requiring horizontal tangents at those points, that is, with zero vertical velocity.
\item \textbf{Cyclic structure:} the pattern forms a smooth closed loop.
\end{enumerate}
The design principles outlined above suggest a separation between the geometric shape of the conducting pattern and the way it is traversed in time.
Our model is built in two stages.

First, we define a single-cycle continuous and smooth geometric curve
\[
g_0:[0,2N]\to\mathbb{R}^2.
\]
where $N$ is the number of beats per cycle.
Because its two endpoint positions and tangents will agree,
it admits a periodic $C^1$ extension
\[
g:\mathbb{R}\to\mathbb{R}^2
\]
with period~$2N$, that is,
\[g(s+2N) = g(s).\]
It describes the geometric trajectory of the baton tip
as a function of a curve parameter~$s$.

Second, we define a single-cycle continuous and smooth timing map
\[
f_0:[0,T]\to[0,2N],
\]
where $T$ is the time per cycle,
and $f_0(0) = 0$ and $f_0(T) = 2N$,
so that $N$~beats take time~$T$.
It is then extended to a global phase function
\[
f:[0,\infty)\to[0,\infty)
\]
by phase accumulation:
\[
f(t)=2N i + f_0\!\left(t-T i\right),\quad \text{where } i = \left\lfloor \frac{t}{T}\right\rfloor.
\]
It maps time~$t$ to the curve parameter~$s$.
The function $f$ satisfies the phase accumulation property
\[
f(t+T)=f(t)+2N,
\]
meaning that each temporal period advances the curve parameter by one full geometric cycle.

The resulting baton motion is
\[
\gamma:[0,\infty)\to\mathbb{R}^2,
\qquad
\gamma(t)=g(f(t)).
\]
Since $g$ has period~$2N$ and $f$ advances by~$2N$ in each temporal period~$T$, it follows that $\gamma$ is periodic with period~$T$.
So, the periodicity of the motion relies on the alignment (synchronization) between the geometric period of $g$ and the phase increment of $f$.
As a side remark,
if these are not commensurate, the resulting trajectory need not be periodic.
In particular, if the phase increment is an irrational multiple of the geometric period, the motion becomes quasi-periodic and does not repeat.

\subsection{Geometric Model}

An $N$-beat pattern is defined by $2N$ \emph{anchor points}
\[P_1, I_1, P_2, I_2, \dots, P_N, I_N,\]
where $P_k$ is the \emph{preparation point} and $I_k$ the \emph{ictus} of beat~$k$.
The sequence is cyclic.
\begin{figure}[hbt]
\centering
{\includegraphics[trim=40mm 30mm 45mm 35mm,clip,height=6cm]{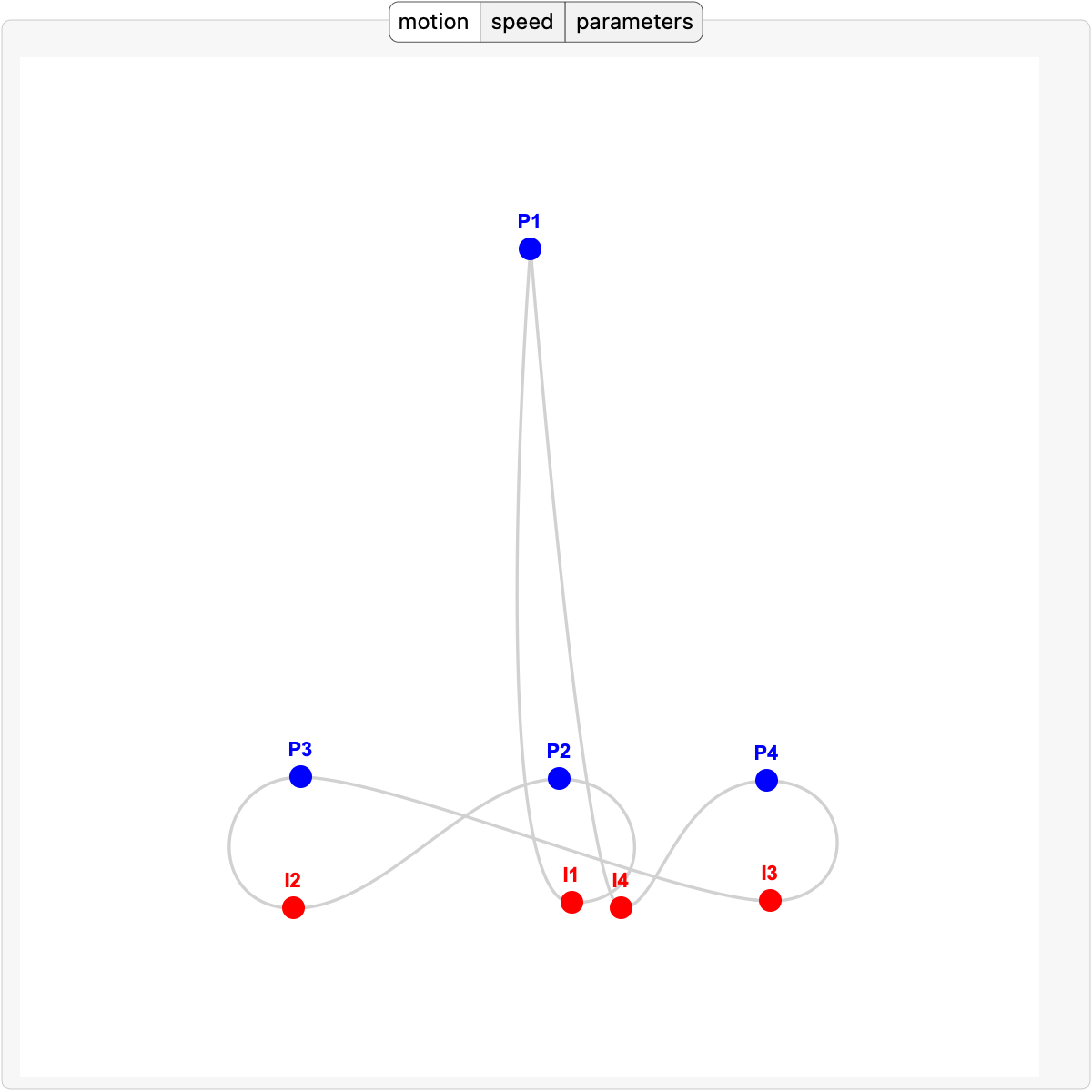}}
\caption{Pattern curve with anchor points at local extremes for 4-beat time signature.}
\label{fig:curve+anchor-points}
\end{figure}

Each anchor point $A_i$ ($i=1,\ldots,2N$,
where $P_k$ corresponds to~$A_{2k-1}$ and $I_k$ to~$A_{2k}$) 
has an associated tangent vector
\[v_i = (r_i, 0),\]
where $r_i \in \mathbb{R}$ is a local signed \emph{roundness parameter}.
Thus, all these tangents are horizontal, ensuring local $y$-extrema at anchors
(see Figure~\ref{fig:curve+anchor-points}).
In general, smaller~$|r_i|$ make the pattern curve sharper (more staccato) and
larger values make it rounder (more legato).
The sign of~$r_i$ determines in which direction the curve passes the anchor point: for $r_i>0$ it goes from left to right
(see points~$I_1, P_3, I_3$ in Figure~\ref{fig:curve+anchor-points}) and
for $r_i < 0$ from right to left
(see points~$P_2, I_2, P_4, I_4$ in Figure~\ref{fig:curve+anchor-points}).
With $r_i=0$ you get a cusp
(see point~$P_1$ in Figure~\ref{fig:curve+anchor-points}).

Between (cyclicly) adjacent anchor points,
the pattern curve is interpolated via cubic Hermite segments~$H_i: [0,1] \to \mathbb{R}^2$ ($i = 1,\ldots,2N$).
In general, for points~$p_0, p_1$ with tangent vectors~$m_0, m_1$,
the curve~$H$ is given by
\begin{align*}
H(u) = {} & (2u^3 - 3u^2 + 1)\,p_0 +
            (u^3 - 2u^2 + u)\,m_0 + {} \\
          & (-2u^3 + 3u^2)\,p_1 +
            (u^3 - u^2)\,m_1, & u \in [0,1].
\end{align*}
This Hermite curve satisfies the following endpoint conditions:
\[
H(0) = p_0, \quad H(1) = p_1, \quad
H'(0) = m_0, \quad H'(1) = m_1.
\]
A cubic is the lowest-degree polynomial that provides the desired modeling freedom.

We take~$H_i$ as the Hermite curve for points~$A_i,A_{i+1}$
and tangents~$v_i,v_{i+1}$ (where $i+1$ is taken cyclicly: $2N+1=1$).
As a consequence, adjacent curve segments meet with matching position and tangent, so the resulting global curve is continuously differentiable
($C^1$; but $C^2$ is not guaranteed), and it is closed.

The full trajectory~$g_0$ of a cycle consists of $2N$~segments~$H_i$
forming a closed curve.
Let $s \in [0,2N]$ be the global curve parameter.
The integer part of~$s$ selects a segment,
and its fractional part gives the local parameter~$u$
for the Hermite cubic:
\[g_0(s) = H_{i+1}(u),
\quad\text{where } i = \lfloor s\rfloor\text{ and }u = s - i.
\]
We typically design a conducting pattern curve
from the conductor's perspective.
The performer’s view is obtained by reflection:
$(x,y) \mapsto (-x,y)$.
To preserve curve shape, roundness must also change sign.

\subsection{Temporal Model}

A full cycle has duration
$T = \frac{60N}{\mathrm{BPM}}$ in seconds.
It is divided into $2N$ equal segments of length
$\Delta = \frac{T}{2N}$.
Odd segments correspond to preparation $\to$ ictus (acceleration),
even segments to ictus $\to$ preparation (deceleration).

Within each segment, time $\tau \in [0,1]$ is mapped to curve parameter $u \in [0,1]$ via a quintic polynomial serving as an ease function:
\[
  \mathrm{ease}(a,b,\tau) = a\tau + c_3 \tau^3 + c_4 \tau^4 + c_5 \tau^5,
\]
with coefficients
\begin{align*}
c_3 & {} = -6a - 4b + 10,\\
c_4 & {} = 8a + 7b - 15,\\
c_5 & {} = -3a - 3b + 6.
\end{align*}
This satisfies:
\begin{align*}
\mathrm{ease}(0) & {} = 0, & \mathrm{ease}(1) & {} = 1,\\
\mathrm{ease}'(0) & {} = a, & \mathrm{ease}'(1) & {} = b,\\
\mathrm{ease}''(0) & {} = 0 \text{ (for free)}, & \mathrm{ease}''(1) & {} = 0.
\end{align*}
A quintic is the lowest-degree polynomial that provides the desired modeling freedom.

\begin{figure}[hbt]
\centering
{\includegraphics[height=5cm]{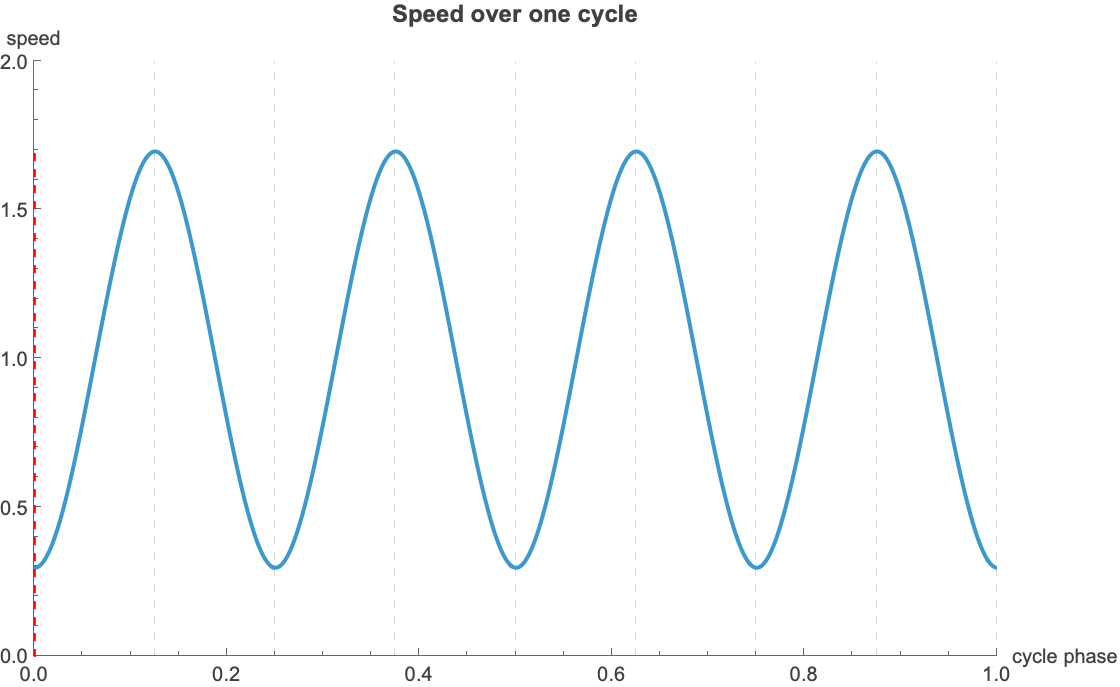}}
\caption{Speed plot for 4-beat pattern with $\beta = 0.7$.}
\label{fig:speed-plot}
\end{figure}
To give us some way to vary the speed along the curve,
we introduce the \emph{speed balance} $\beta \in [0,1]$.
Defining
\[a_{\min} = 1 - \beta,\quad a_{\max} = 1 + \beta,\]
we then have
\begin{itemize}
\item prep $\to$ ictus: $a = a_{\min},\ b = a_{\max}$,
\item ictus $\to$ prep: $a = a_{\max},\ b = a_{\min}$.
\end{itemize}
with interpretation
\begin{itemize}
\item $\beta = 0$: uniform motion,
\item $\beta = 1$: maximal contrast (zero at preparation, fastest at ictus).
\end{itemize}
Figure~\ref{fig:speed-plot} shows the speed plot for a 4-beat pattern
with $\beta=0.7$.

Let $t \in [0,T]$ be the global time parameter.
The integer part of~$t/\Delta$ selects a segment,
and its fractional part gives the local parameter~$\tau$
for the quintic ease function:
\[f_0(t) = \mathrm{ease}_{i+1}(u),
\quad\text{where } i = \lfloor t/\Delta\rfloor\text{ and }\tau = t/\Delta - i.
\]

\subsection{Implementation}

\begin{figure}[hbt]
\centering
\includegraphics[width=\textwidth]{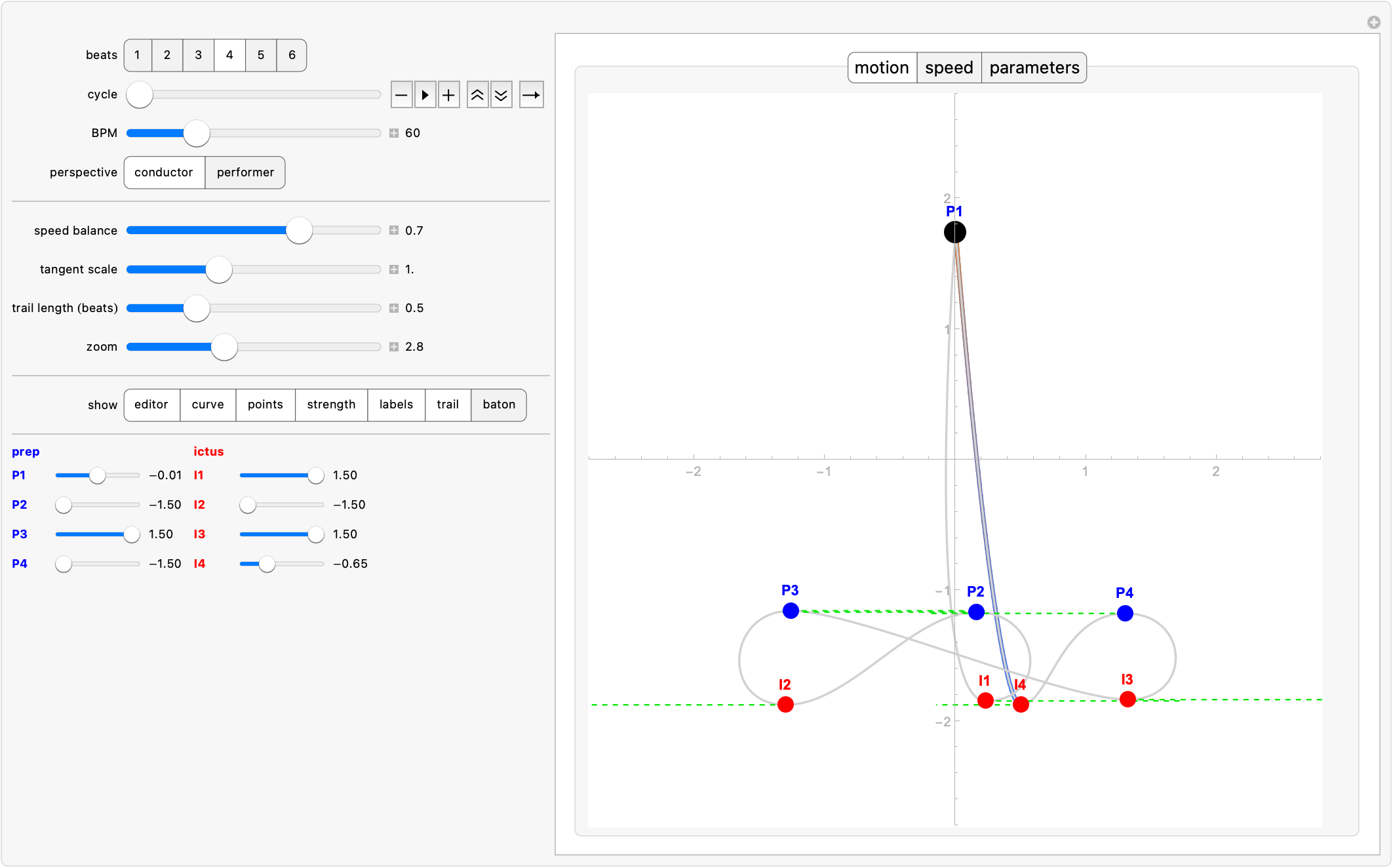}
\caption{Screenshot of our Wolfram Demonstration \emph{Conducting Patterns}.}
\label{fig:demo}
\end{figure}

The model has been implemented as an interactive Wolfram Demonstration \emph{Conducting Patterns}~\cite{Verhoeff:Wolfram}
(also see Figure~\ref{fig:demo}), allowing:
dragging of anchor points,
adjustment of roundness (tangent magnitude),
control of tempo and speed balance, and
visualization of trajectory and time map.
This demonstration was used to create the illustrations of this article.

The model also features in the web app \emph{Crusis},
a conducting companion~\cite{Verhoeff:Crusis}.
It offers a pattern editor based on this model.
\emph{Crusis} incorporates some enhancements
beyond a dot traversing the pattern curve,
such as
\begin{itemize}
\item drawing a fat dot at the baton's tip location on the pattern curve,
\item changing the color of that fat dot momentarily at the downbeat,
\item drawing the baton as a line from a fixed origin
  to the point on the pattern curve, and 
\item drawing a trail behind the baton tip along the pattern curve.
\end{itemize}

\section{Related Work}

\subsection{Geometric modeling of trajectories}

The geometric component of our model draws on classical results from curve design. In particular, cubic Hermite interpolation provides a local, endpoint-constrained representation with explicit control over tangent directions and magnitudes \cite{Farin,Prautzsch}. This makes it especially suitable for our setting, where horizontal tangents at anchor points (preparations and icti) are a structural requirement. Compared to global spline constructions (e.g., B-splines), the Hermite form offers strictly local control, avoiding unintended global deformations when adjusting a single beat.

Local interpolating splines, such as Catmull--Rom splines \cite{CatmullRom,Farin},
share the property of locality but typically infer tangents from neighboring points.
Our model instead prescribes tangents explicitly,
reflecting the conducting constraint that extrema must have horizontal tangents.
This explicit control is essential for guaranteeing that preparations and icti are visually salient as local maxima and minima.

\subsection{Temporal parametrization and easing}

The temporal component uses a quintic polynomial to map normalized time within each half-beat to curve parameter. Quintic ``ease'' functions are standard in animation for achieving $C^2$-smooth motion with prescribed endpoint velocities, thereby avoiding visible kinks in velocity and acceleration \cite{Farin}. In our context, the endpoint derivatives correspond to minimal and maximal phase speeds, allowing a single parameter to control the contrast between preparation and ictus.

More broadly, time reparameterization is a standard idea in control and motion planning, where one distinguishes a geometric path from the law according to which that path is traversed \cite{Bobrow+al}. Our approach adapts this idea to conducting: the geometry defines the path, while a low-dimensional timing law modulates how the baton traverses it. This separation enables independent reasoning about spatial form and temporal emphasis.

\subsection{Empirical and motion-capture studies}

Empirical work on conducting gestures has focused on measurement and analysis rather than minimal modeling. Marrin Nakra \cite{Nakra} developed sensor-based systems to capture conductor motion, revealing characteristic patterns of acceleration and emphasis.
Luck and Toiviainen \cite{LuckToiviainen} studied synchronization between conductors and ensembles, highlighting the perceptual importance of the ictus. 
Wöllner \cite{Wollner} surveys gesture in music performance more broadly, emphasizing its communicative and expressive roles.

These studies support our design principles: the prominence of ictus points, the preparatory motion leading into them, and the importance of smooth yet expressive timing. However, they do not provide a compact parametric model.
The present work can be seen as a complementary, structural counterpart to these empirical approaches.

The scientific study of conducting gesture is relatively recent and remains less developed than motion analysis in domains such as sports or surgery.
Motion-capture studies,
such as the ``Capturing the Contemporary Conductor'' project~\cite{Polfreman+al},
have used three-dimensional tracking to record detailed hand, body, and facial movements of conductors.
These studies aim to provide high-resolution data on conducting technique,
but they do not typically propose compact mathematical models.

A consistent empirical finding across such studies is that the beat (ictus) corresponds to a local minimum in the vertical trajectory of the baton tip. That is, the beat is visually identified at a point where the motion reverses direction. This observation directly supports the modeling constraint used in the present work, where ictus points are represented as local minima of the trajectory with horizontal tangents.

From a biomechanical perspective, the change in direction at the ictus may be produced using different joints (shoulder, elbow, or wrist), with the elbow often serving as the primary hinge for basic timekeeping. The baton itself functions as an extension of the arm, with the ictus ideally located at the tip of the stick.

Conversely,
the preparation point---the high point of the rebound arc preceding each ictus---corresponds to a local maximum in the vertical trajectory,
as described in standard conducting pedagogy~\cite{Atherton, Rudolf}.

\subsection{Gesture recognition}

A substantial body of work in computer science has focused on recognizing and following conducting gestures using sensors.
Systems have been developed using devices such as the Nintendo Wii Remote, Microsoft Kinect, and custom digital batons equipped with accelerometers.
These systems aim to detect tempo, identify beat positions, and classify expressive features such as articulation~\cite{Toh+al}.

Although primarily engineering-driven, these approaches implicitly define geometric and temporal features of conducting gestures. For example, beat detection algorithms often rely on identifying extrema or direction changes in motion trajectories, aligning with the empirical characterization of the ictus as a local minimum.

However, these systems typically model gestures indirectly, through feature extraction and classification, rather than providing an explicit parametric description of the underlying motion. The present work complements this line of research by proposing a direct mathematical model of conducting patterns.

\subsection{Conducting pedagogy}

In conducting pedagogy, standard beat patterns are taught as canonical shapes, with emphasis on clarity of ictus and preparation~\cite{Rudolf}.
The ictus is defined as the point of direction change that communicates the beat to the ensemble, while the preparation guides performers toward it.

Pedagogical sources emphasize that the baton tip should clearly indicate the ictus, and that effective conducting involves a balance between clarity and fluidity of motion. Gesture height, conducting plane, and the choice of joints used in the motion all influence how performers perceive timing and expressivity.

These pedagogical principles are reflected in the constraints of the present model, particularly the requirement of horizontal tangents at anchor points and the separation between geometric trajectory and temporal emphasis.

\section{Discussion}

The presented model integrates ideas from geometric modeling,
temporal reparameterization, and empirical observations of conducting.
The model achieves:

\begin{itemize}
\item a strict separation between geometry and timing,
\item a minimal parameter set
  (for an $N$-beat pattern, $6N$~real numbers for the geometry,
  and speed balance~$\beta\in[0,1]$ for timing),
\item smooth motion with controlled extrema,
\item expressive flexibility via the single timing parameter~$\beta$.
\end{itemize}
The separation between geometry and timing clarifies the roles of different parameters: anchor points and roundness determine the visual shape, while the timing law determines how that shape is traversed. This mirrors conducting practice, where the shape of the gesture and its temporal emphasis are conceptually distinct.

To our knowledge, there is no existing model that explicitly separates geometry and timing for conducting patterns in this way, with both components expressed in closed form and controlled by a small number of interpretable parameters.

The model is purely geometric and kinematic: it does not derive from an underlying physical or biomechanical model of the arm, and makes no claims about the forces or joint movements involved in real conducting.

Certain parameter choices may produce undesirable artifacts:
\begin{itemize}
\item the geometric curve is only~$C^1$ at anchor points;
  $C^2$ continuity---and thus smooth curvature---would require either
  higher-degree segments or additional constraints on the tangent magnitudes.
\item non-uniform spatial speed due to parameterization.
\end{itemize}
It is important to note that the temporal model governs the parameterization of the curve rather than its arc length. Consequently, uniform phase speed does not imply uniform physical velocity of the baton tip.
This distinction explains certain perceptual artifacts and
suggests arc-length reparameterization as a possible refinement.

\section{Conclusion}

We have presented a compact mathematical model for conducting patterns
that captures essential geometric and temporal features while remaining computationally simple.
The model builds on the well-established observation that
the ictus occurs at a local vertical minimum, and that the baton tip accelerates toward the ictus and decelerates away from it.
We add the explicit modeling of the preparation point as a local vertical maximum,
paired with the ictus at equal temporal distance,
forming the two anchor points of each beat.
The model is purely geometric and kinematic:
it makes no claims about the underlying physics or biomechanics of the conducting arm.
What appears to be novel is the explicit formal requirement that
the preparation falls at the exact temporal midpoint between consecutive icti.

Future work includes:
\begin{itemize}
\item arc-length reparameterization,
\item extension to three-dimensional motion,
\item empirical validation with musicians.
\end{itemize}

\appendix

\section{Example 2-, 3-, and 6-beat conducting patterns}

To illustrate the use of our model,
we show, in addition to the 4-beat pattern in Figure~\ref{fig:curve+anchor-points},
parameterizations of the common 2-, 3-, and 6-beat time signatures.
The choice of anchor point positions and roundness values
used as defaults in the implementation
was informed by conducting practice,
as illustrated in instructional videos by Baumeister~\cite{Baumeister}.

In general, a timing balance~$\beta$ between~$0.5$ and~$0.7$ works well.
\begin{figure}[hbt]
{\includegraphics[trim=80mm 50mm 60mm 35mm,clip,height=6cm]{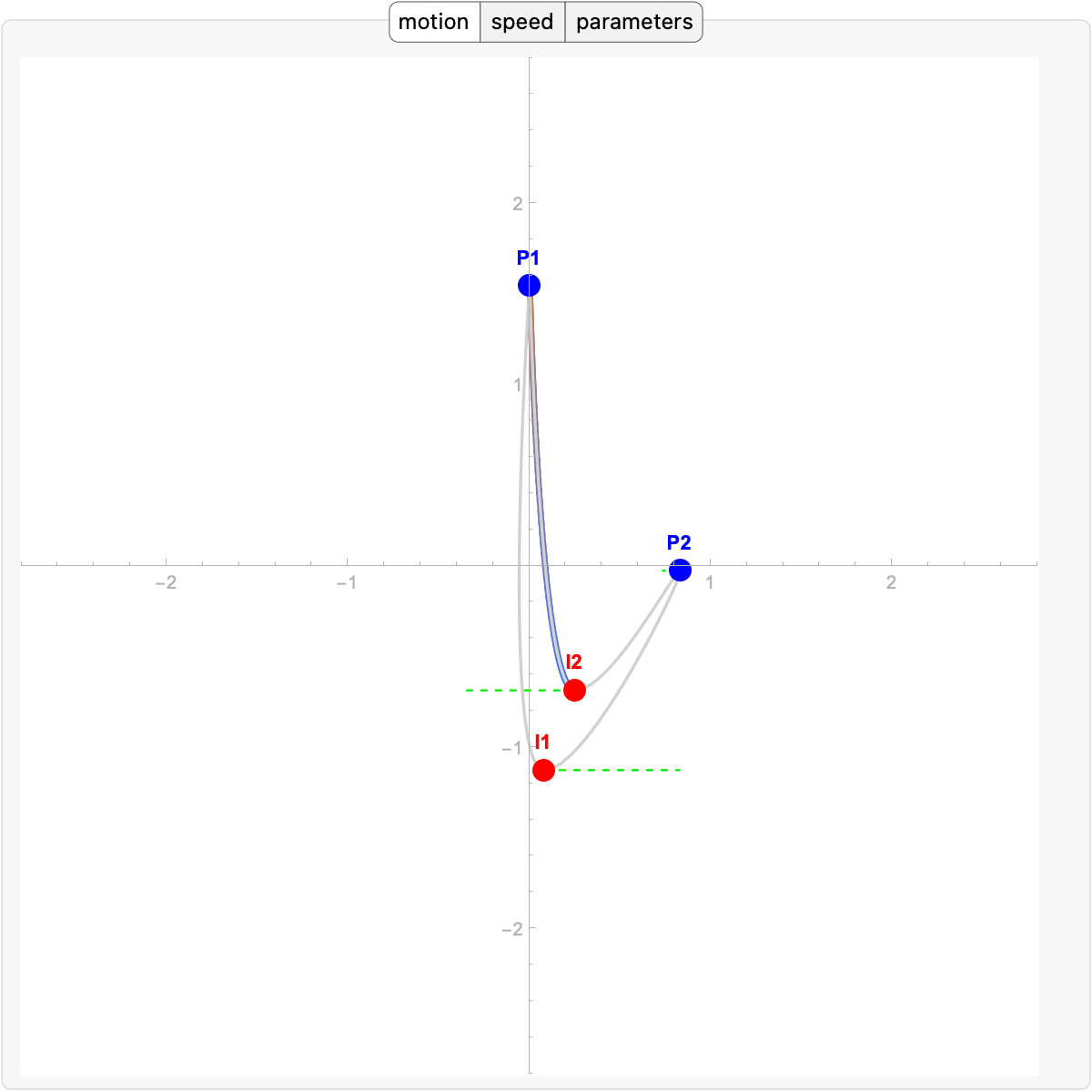}}\hfil
{\includegraphics[trim=60mm 30mm 45mm 35mm,clip,height=6cm]{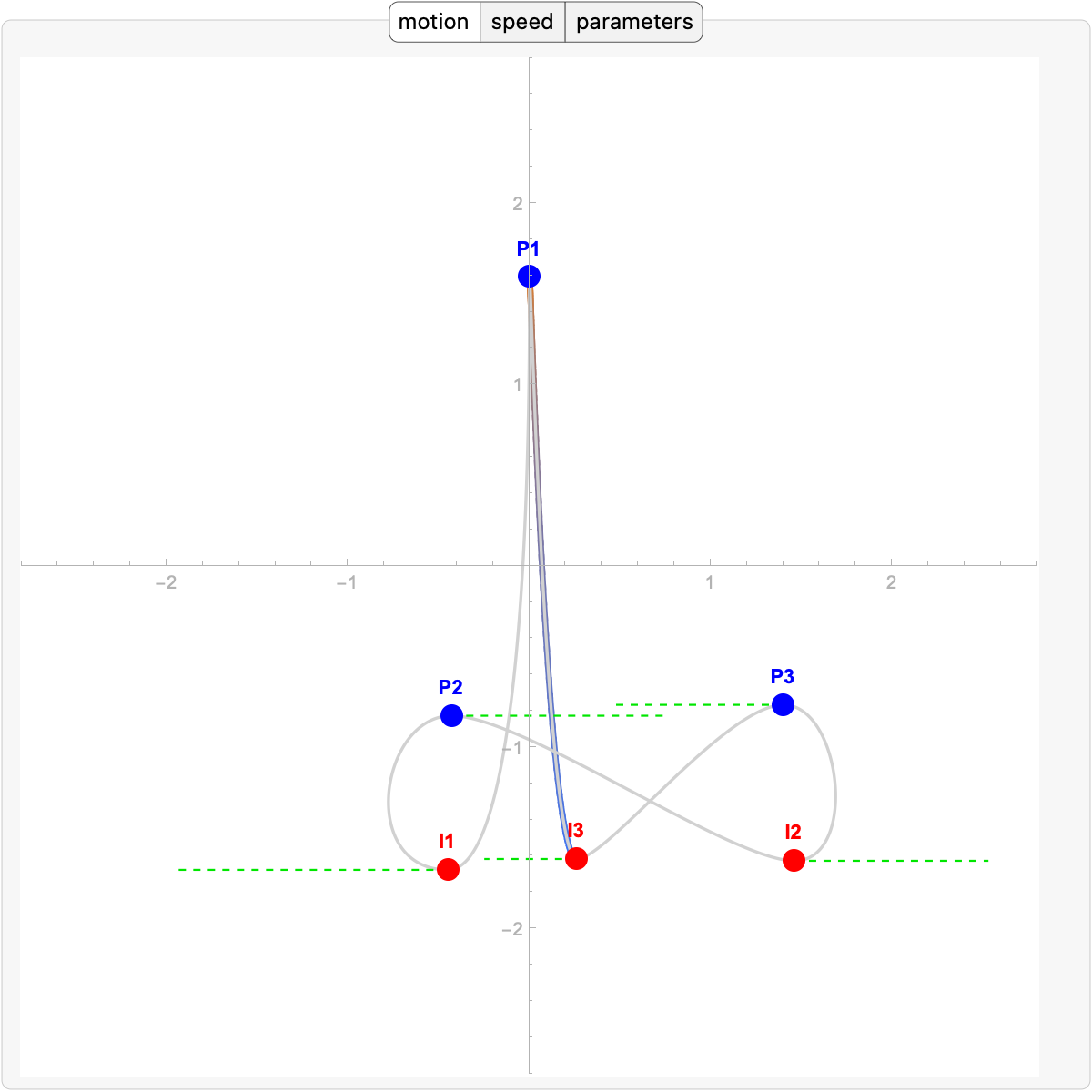}}\hfil
{\includegraphics[trim=40mm 20mm 40mm 35mm,clip,height=6cm]{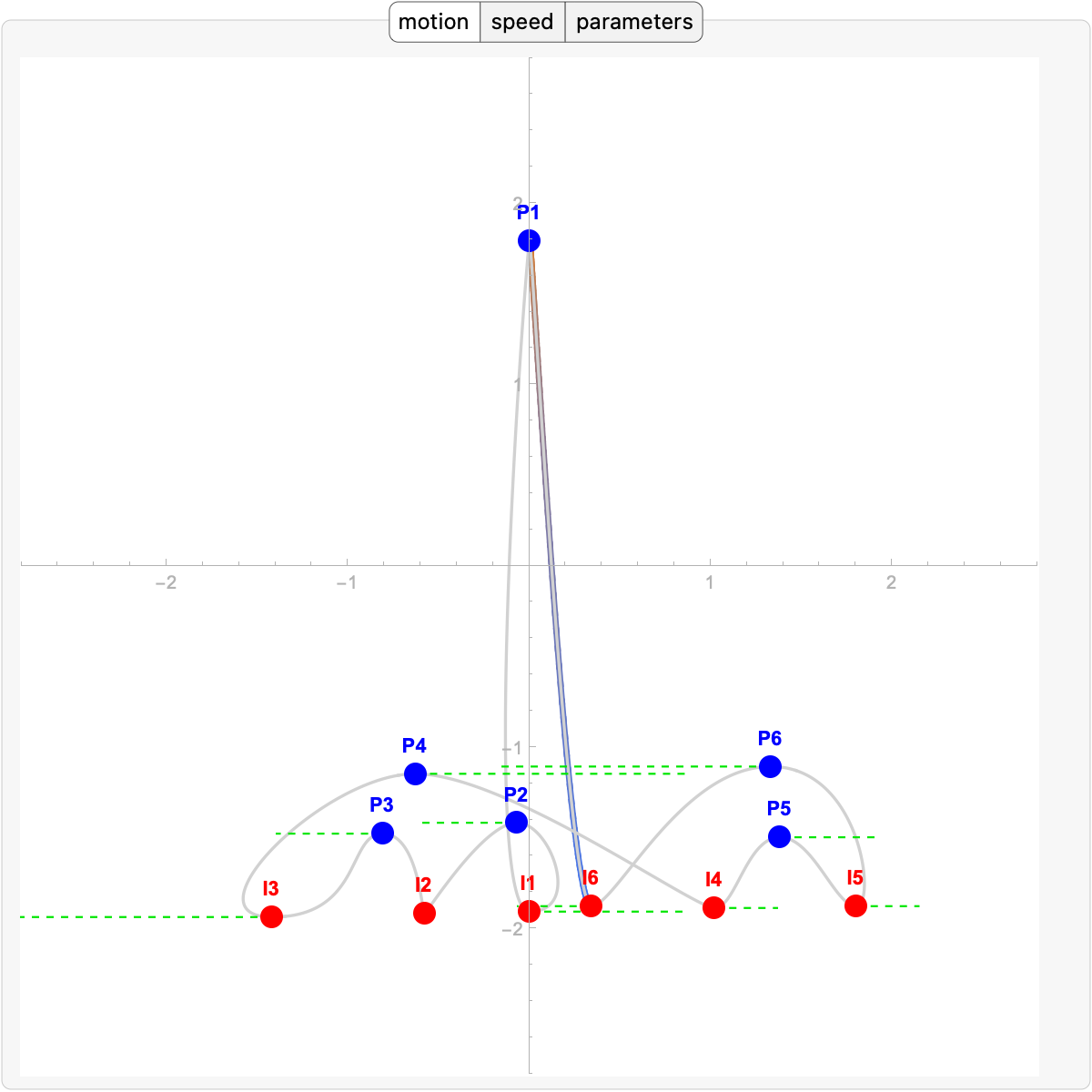}}
\caption{Common pattern curves for 2-, 3-, and 6-beat time signatures.}
\label{fig:2-3-6-beat-pattern-curves}
\end{figure}

\end{document}